\documentclass[11pt, twoside, fleqno]{article} %%, draft

\usepackage{stmaryrd}
\usepackage{amsmath}
\usepackage{a4}
\usepackage[english, francais]{babel}
\usepackage{euscript}
\usepackage{amsfonts}
\usepackage{amssymb}
\usepackage{mathrsfs}
\usepackage{xypic}
\usepackage[dvips]{graphics}
\usepackage{enumerate}
\usepackage{array}
\usepackage{theorem}
\xyoption{all}
\usepackage[T1]{fontenc}
\usepackage{makeidx}
\usepackage{float}

\newcommand{\Spec}{{\operatorname{Spec}\kern 1pt}}
\newcommand{\Hom}{{\mathrm{Hom}}}
\newcommand{\p}{{\mathfrak p}}
\newcommand{\qed}{\hfill $\Box$ \medskip}
\newcommand{\preuve}{\noindent{\it Preuve : }}
\newcommand{\GG}{\ensuremath{\mathbb G}}
\newcommand{\PP}{{\mathbb P}}
\newcommand{\ZZ}{{\mathbb Z}}
\newcommand{\Aut}{{\mathrm{Aut}}}
\newcommand{\Ext}{{\mathrm{Ext}}}
\newcommand{\Frac}{\mathrm{Frac}}
\renewcommand{\H}{{\mathrm{H}}}
\newcommand{\aalpha}{\ensuremath{\boldsymbol \alpha}}
\newcommand{\dif}{{\mathfrak d}}

\theoremstyle{plain}
\newtheorem{theorem}{Théorème}[section]
\newtheorem{proposition}[theorem]{Proposition}
\newtheorem{corollary}[theorem]{Corollaire}
\newtheorem{lemma}[theorem]{Lemme}

\theorembodyfont{\upshape}

\newtheorem{definition}[theorem]{Définition}

\newtheorem{remark}[theorem]{Remarque}

\begin{document}

\title{Aller-retour vers l'inséparable}

\author{Sylvain Maugeais}

\date{}

\maketitle

\begin{abstract}
Nous construisons des morphismes inséparables entre courbes de genre
$\ge 2$ qui proviennent par dégénérescence de morphismes séparables.
\end{abstract}

Soient $R$ un anneau de valuation discrète de caractéristique
résiduelle $> 0$ et $f\colon \mathcal X \to \mathcal Y$ un
$R$-morphisme fini génériquement séparable entre courbes propres et
lisses. Se peut il que le morphisme induit par $f$ au-dessus du point
fermé de $\Spec R$ soit inséparable ? En d'autres termes, peut on
dégénérer du séparable vers l'inséparable ? Lorsque le genre de
$\mathcal X$ est $\le 1$, il est aisé de construire des exemples mais
lorsque ce genre est $\ge 2$ cela devient plus compliqué. Par exemple,
Deligne et Mumford ont montré dans \cite{Deligne_Mumford} que cela ne
pouvait pas se passer dans le cas d'un morphisme $f$ galoisien.

Le but de cet article est de donner un critère simple (dans le cas
où le degré d'inséparabilité est petit) pour répondre à cette
question, puis de produire des exemples explicites.

L'étude classique des déformations de morphismes ne donne
malheureusement pas beaucoup d'information à ce sujet. En effet, les
obstructions à la déformation ne sont pas locales dans le cas d'un 
morphisme inséparable, ceci rend l'étude particulièrement difficile.

Pour contourner ce problème, nous introduisons dans le présent article
la notion de morphismes rigidifiés qui consiste en la donnée d'un
morphisme $f$ entre courbes lisses et d'une décomposition de
l'homomorphisme $df$, induit sur les faisceaux de différentielles,
sous la forme $df=\xi\delta$ où $\xi$ est une constante (éventuellement
nulle) et $\delta$ est un homomorphisme \emph{injectif dans chaque
fibre} et qui envoie les formes exactes sur des formes exactes. Cette
condition technique intervient naturellement dans l'étude des
dégénérescences de morphismes entre courbes lisses en \emph{égales
  caractéristiques} et permet ensuite de répondre au cas de l'inégale
caractéristique par un argument global sur des espaces de modules.

Une fois cette définition posée, nous procédons à l'étude des
déformations des morphismes rigidifiés ayant un degré d'inséparabilité
au plus égal à la caractéristique et montrons que ce problème est
non obstrué (cf. Théorème \ref{DeformationMorphismeRigidifie}). Un
corollaire immédiat nous permet de démontrer le théorème suivant 

\begin{theorem}
Soit $f\colon X \to Y$ un morphisme fini entre courbes propres et
lisses au-dessus d'un corps de caractéristique $p$. Supposons que $f$
se décompose en $h \circ F$ où $F$ est purement inséparable de degré
$p$ et $h$ est séparable. 

Alors $f$ provient par dégénérescence en égale caractéristique d'un
morphisme séparable si et seulement si $\Hom(h^* \Omega_{Y/k}, B^1_X)
\not = \{0\}$ ($B^1_X$ désignant le faisceau des différentielles
exactes). 
\end{theorem}

En particulier, on peut constater que ce théorème ne dit rien sur les
morphismes dont le degré d'inséparabilité est $> p$, non plus que sur
les déformations vers la caractéristique $0$. Ce deuxième problème
peut être contourner utilisant un raisonnement sur des espaces
de modules, mais on obtient alors seulement une condition nécessaire 
(cf. corollaire \ref{RelevementCaracteristitique0}). Par contre, nous
ne pouvons pas, à l'heure actuelle, déformer des morphismes de gros
degré d'inséparabilité.

Une fois ce critère en main, le problème est de montrer qu'il existe
des morphismes de degré d'inséparabilité $p$ qui admette une structure
de morphisme rigidifié, ce qui se trouve être une condition très
forte. Il est assez aisé de trouver des conditions nécessaire
(cf. section \ref{CN}) mais plus difficile de trouver des conditions
suffisantes, ces ce qui est fait dans la dernière partie de ce
travail, cf. exemple \ref{CS}

\noindent{\bf Remerciements :} Je souhaite remercier Qing Liu pour
m'avoir posé la question des relèvements de morphismes inséparables.

%\subsection{Préliminaires}
%
%\begin{proposition}
%Soit $f\colon X \to Y$ un morphisme fini entre courbes propres et
%lisses sur un corps de caractéristique $p$. Alors $f$ se décompose en
%$X \to X^{(p^r)}\to Y$ où $X \to X^{(p^r)}$ est une puissance du
%Frobenius relatif et $X^{(p^r)} \to Y$ est séparable.
%
%De plus, si $k$ est parfait alors $f$ se décompose en 
%$X \to Z \to Y$ avec $X \to Z$ séparable et $Z \to Y$ une puissance du
%Verschiebung relatif.
%\end{proposition}
%
%\preuve D'après Lang, on peut décomposer en purement inséparable +
%séparable. Le morphisme purement inséparable est alors une puissance
%du Frobenius  d'après \cite{Liu}, Proposition 7.4.21. 
%
%Considérons $Z$ la courbe propre et lisse dont le corps des
%fractions est $K(Y)^{1/p^r}$ qui existe car $k$ est parfait 
%(et donc le Frobenius possède un inverse, en effet, si $Y$ est donnée
% localement par $k[x_1, \ldots, x_n]/I$ alors $Y^{1/p^r}$ est donnée
% par $k[x_1, \ldots, x_n]/I'$ avec $I'$ engendré par les $\sum
% a_i^{1/p^r} x^i$ avec $\sum a_i x^i \in I$, $a_i^{1/p^r}$ n'a un sens
%en général que si $k$ est parfait). 
%En particulier, $K(X)=K(X^{(p^r)})^{1/p^r}$ et donc (via $f$) on a $K(Z) \subset K(X)$.
%
%On a une tour d'extension 
%$$\xymatrix{
%& K(X^{(p^r)}) K(Z) \subset K(X)\ar[rd] \ar[ld] & \\
%K(X^{p^r}) \ar[rd] & & K(Z) \ar[ld] \\
%& K(Y) &
%}$$
%on conclut alors en regardant les degrés. \qed

\section{Définition}

Dans tout cet article, les démonstrations se font souvent en se
ramenant au cas des séries formelles. Ceci est justifié par le lemme
suivant.

\begin{lemma}
Soit $S$ un schéma d'égale caractéristique $p > 0$ et $X \to S$ une courbe
lisse. Soit $\omega$ une forme différentielle sur $X$. Alors 
\begin{enumerate}[i)]
\item $\omega$ est localement exacte si et seulement si elle l'est sur
  $X \times_S S'$ pour $S' \to S$ fidèlement plat.
\item $\omega$ est localement exacte si et seulement si pour tout
  point $\p \in X$, l'image de $\omega$ dans
  $\Omega_{X/S}\otimes_{\O_X} \widehat O_{X, \p}$ est exacte.
\end{enumerate}
\end{lemma}

Pour montrer ce lemme, il convient d'introduire quelques notations : 
notons  $F\colon X \to X^{(p)}$ le Frobenius relatif et $B^1_{X/S}$
l'image de $d\colon \O_X \to \Omega_{X/S}$. Alors 
le faisceau $F_* B^1_{X/S}$ est naturellement un sous-faisceau
cohérent plat de $F_* \Omega_{X/S}$ : c'est le noyau de l'opérateur de
Cartier relatif $C \colon F_* \Omega_{X} \to \Omega_{X^{(p)}}$. On
identifiera souvent $B^1_{X/S}$ a son image via $F_*$.

Le lemme ci-dessus est alors une conséquence immédiate de la cohérence
de $F_* B^1_{X/S}$.

\begin{definition}
Soit $f\colon X \to Y$ un morphisme fini entre $S$ courbes lisses (non 
nécessairement propres). Un homomorphisme $\delta\colon 
\Omega_{Y/S} \to f_*\Omega_{X/S}$ sera dit $f$-exact s'il induit un
homomorphisme $B^1_{Y/S} \to f_* B^1_{X/S}$ (i.e. il envoie une forme
localement exacte sur une forme localement exacte).
\end{definition}

L'exemple le plus simple d'homomorphisme exact est l'homomorphisme
$df$. Pour toute forme localement exacte $\omega=dh$ on 
a $df(\omega)=d(f(h))$.

\begin{lemma}
\label{FaisceauHomfExct}
Supposons que $S$ est le spectre d'un corps et que $f$ est
inséparable. Alors l'ensemble des homomorphismes $f$-exact s'identifie
à $\Hom(\Omega_Y, f_*B^1_X)$. 

Si on a de plus une décomposition $f=h \circ F$, alors 
$$\Hom(\Omega_Y, f_*B^1_X)=\Hom(h^* \Omega_Y, F_*B^1_X).$$
\end{lemma}

\preuve Soit $\delta$ un homomorphisme $f$-exact. Au voisinage de
chaque point de $Y$, il existe une base $\omega$ de $\Omega_Y$ qui est  
localement exacte. En particulier, $\delta(\omega)=dx \in f_* B^1_{X/S}$.
Comme $f$ est inséparable, pour tout élément $\alpha \in \O_Y$ on a
$d(f(\alpha))=0$ et donc $\delta(\alpha \omega)=f(\alpha)
\delta(\omega)=f(\alpha) dx=d(f(\alpha) x)$. Donc $\delta$ est à
valeurs dans $f_* B^1_{X/S}$.

La dernière assertion vient par adjonction. \qed

\begin{definition}
Un morphisme rigidifié au-dessus d'un anneau local $A$ consiste en la donnée 
d'un triplet $(f, \delta, \xi)$ où $f$ est un morphisme fini entre
$A$-courbes propres et lisses, $\delta$ un homomorphisme $f$-exact non
nul au-dessus de chaque point de $A$, $\xi\in A$ tels que $df=\xi\delta$.
\end{definition}

Il existe une action naturelle de $\GG_m$ sur les morphismes rigidifiés
donnée par $\alpha.(f, \delta, \xi)=(f, \alpha \delta,
\alpha^{-1}\xi)$, ce qui complique une éventuelle définition de ce que
pourrait être des morphismes rigidifiés sur une base
globale (il suffit toutefois de considérer le champ quotient afin
d'obtenir un objet global). Nous n'utiliserons ici que des considérations
locales ! 

\medskip

Le fait que $\delta$ soit non nul en chaque point équivaut en fait à ce que 
l'homomorphisme $f^* \Omega_Y \to \Omega_X$ induit par adjonction ait un conoyau 
fini plat sur $A$. En particulier la notion de morphisme rigidifié est
stable par changement de base.

\begin{proposition}
Le morphisme sous-jacent d'un morphisme rigidifié $(f, \delta, \xi)$
est séparable au-dessus d'un point $\p \in \Spec A$ si et seulement si
$\xi$ est inversible en $\p$. 
\end{proposition}

Un morphisme $f$ entre $A$-courbes propres et lisses qui est séparable
au-dessus de chaque point de $A$ possède une structure canonique de
morphisme rigidifié : $(f, df, 1)$. 

\section{Dégénérescence}

La notion de morphisme rigidifié est justifiée par la proposition
suivante.

\begin{proposition}
Soit $R$ un anneau de valuation discrète d'\emph{égales
caractéristiques $p$} et $f\colon X \to Y$ un morphisme fini entre 
$R$-courbes propres et lisses. Supposons que $f$ est séparable à la
fibre générique.  
Alors $f$ peut être muni d'une structure de morphisme rigidifié.
\end{proposition}

\preuve Considérons la partie verticale du diviseur de ramification :
c'est un multiple de la fibre spéciale de $X$. Comme celle-ci
est réduite, ce diviseur est de la forme $(\xi)$ avec $\xi \in R$. Posons
$\delta=\frac{df}{\xi}$. Comme $R$ est d'égale caractéristique, $\delta$
est $f$-exacte. Pour montrer cette dernière assertion, on se ramène au
cas d'un morphisme $f\colon R[[y]] \to R[[x]]$ par changement de base
fidèlement plat, localisation et complétion. Par suite, comme $R$ est
d'égales caractéristiques, un élément $\omega=(\sum a_i x^i)dx$ est
exact si et seulement si $a_i=0$ pour tout $i=-1 \mod p$. Ceci impose
que $\frac{\omega}{\xi}$ est exacte. \qed

Soit $R$ un anneau de valuation discrète d'inégales caractéristiques
$(0, p)$ de valuation $v$ et $f\colon X \to Y$ un morphisme fini ente
$S$-courbes propres et lisses et notons $p^i$ le degré
d'inséparabilité à la fibre spéciale. Notons $\xi$ un générateur de la
partie verticale du diviseur de ramification. Comme dans le cas
d'égales caractéristiques, l'homomorphisme $\frac{df}{\xi}$ a encore un
sens mais n'est plus exacte en général. Toutefois, $\xi$ ne peut par être
quelconque. En effet, notons $v$ la valuation sur $R$.

\begin{lemma}
On a $v(\xi) \le  v(p^i)$.
\end{lemma}

\preuve Se placer en un point de la fibre spéciale en lequel $f$ se
décompose en étale + purement inséparable et se ramener ainsi à
l'étude de $f:R[[y]] \to R[[x]]$ où $f$ est purement inséparable à la
fibre spéciale. Il suffit alors de regarder $d(f(y))$. \qed

\begin{proposition}
Supposons que le degré d'inséparabilité à la fibre spéciale soit $p$
et que $v(\xi) < v(p)$.  
Alors il existe un homomorphisme $f$-exact non trivial.
\end{proposition}

\preuve Il suffit de montrer que $\frac{df}{\xi}$ est $f$-exact. 
Pour cela, on se ramène, par localisation  puis completion au cas des
séries formelles $R[[x]]$ et on le démontre explicitement. \qed  

Si on retire l'hypothèse que $v(a) < v(p)$, alors ce 
résultat n'est plus vrai, comme le montre l'exemple suivant :
considérons le morphisme $\PP^1_{\ZZ_p} \to \PP^1_{\ZZ_p}$ définit par $x \mapsto
x^p$. Ce morphisme ne donne pas naissance à un homomorphisme exact
à la fibre spéciale !

\section{Théorie de la déformation}

Le but de cette section est de démontrer le théorème suivant.

\begin{theorem}
\label{DeformationMorphismeRigidifie}
Soient $(f_0\colon X_0 \to Y_0, \delta_0, 0)$ un morphisme rigidifié
entre courbes propres et lisses au-dessus d'un corps parfait $k$, $A$ 
un anneau  noethérien, local complet de caractéristique $p$ et de
corps résiduel $k$, et $\xi \in A$. Si le degré d'inséparabilité de
$f_0$ est au plus $p$, alors il existe un morphisme rigidifié
$(f\colon X \to Y, \delta, \xi)$ qui soit une déformation de $(f_0,
\delta_0, 0)$. 
\end{theorem}

Le cas où $f$ est séparable est bien connu : c'est un corollaire assez
direct de la théorie des déformations des morphismes telle qu'exposée
dans \cite{ZivRan}. Dans la suite, nous utiliserons cette théorie
librement (en particulier la classification des automorphismes de
déformations). 

Le cas séparable étant réglé, nous supposerons désormais que le degré
d'inséparabilité de $f$ à la fibre spéciale est exactement $p$. En
particulier, le morphisme est inséparable.

On se donne donc $(f_0\colon X_0 \to Y_0, \delta_0, 0)$ un morphisme
rigidifié au-dessus d'un corps $k$ et $(f\colon X \to Y, \delta, \xi)$
une déformation au-dessus d'un anneau artinien $A$ de caractéristique
$p$ et de corps résiduel $k$.

\begin{lemma}
\label{Automorphisme}
Soit $A' \to A$ une petite extension et $(\tilde f\colon \tilde X \to
\tilde Y, \tilde \delta, \tilde \xi)$ une déformation de $(f\colon X
\to Y, \delta, \xi)$. Alors le groupe des automorphismes de cette
déformation est naturellement isomorphe au noyau du morphisme
$$\begin{array}{ccc}
\Hom(\Omega_{X_0}, \O_{X_0}) & \to & \Hom(\Omega_{Y_0}, \Omega_{X_0}) \\
\chi & \mapsto & d \circ \chi \circ \delta.
\end{array}$$
En particulier, son image directe par le Frobenius absolu est cohérente
(et donc sa cohomologie sur un schéma affine est triviale).
\end{lemma}

\preuve A priori, un automorphisme de déformation de $(f, \delta,
\xi)$ est un élément $(\sigma, \tau) \in \Aut(\tilde X) \times
\Aut(\tilde Y)\cong \Hom(\Omega_{X_0}, \O_{X_0}) \times
\Hom(\Omega_{Y_0}, \O_{Y_0})$. En utilisant la compatibilité $f=\tau
\circ f \circ \sigma$ et le fait que $f$ est inséparable à la fibre
spéciale, on trouve que $\tau = Id$.

Le résultat s'ensuit en écrivant explicitement l'action de
$\Hom(\Omega_{X_0}, \O_{X_0})$ sur $\delta$. \qed

\begin{lemma}
\label{ClasseDeformation}
Supposons $A$ artinien, $X$ et $Y$ affine, $\delta$ est un
isomorphisme et à la fibre spéciale le morphisme $f_0$ est la composée
d'un morphisme étale et d'un morphisme purement inséparable de degré
$p$. Donnons nous une petite extension $A' \to A$ de noyau
$(\epsilon)$ et deux relèvements $(f_1, \delta_1, \xi)$ et $(f_2,
\delta_2, \xi)$ (on a le même $\xi$). Alors il existe des
automorphismes de déformations $\phi \in \Aut_X(X_1)$ et $\psi \in
\Aut_Y(Y_1)$ tels que $f_2=\psi \circ f_1 \circ \phi$ et
$\delta_1=d\phi \circ \delta_2 \circ d\psi$. 
\end{lemma}

\preuve Comme $X$ et $Y$ sont lisses et affines, on peut supposer que
$X_1=X_2$ et $Y_1=Y_2$. 

Choisissons une base locale $dx$ de $\Omega_X$ et une base $dy$ de
$\Omega_{Y_0}$ (ce qui peut toujours se faire, quitte à réduire $X$ et
$Y$) de sorte que $\Ext^0(\Omega_{X_0}, \O_{X_0}) = \O_{X_0}
\frac{\partial}{\partial x}$.  On a
$\delta_1(dy)-\delta_2(dy)=\epsilon Q dx$ avec $Q \in
\O_{X_0}$. Notons $\delta_1(dy)=\alpha dx$. Par hypothèse, $\delta_1$
est un isomorphisme donc $\alpha$ est inversible. Posons
$h_1=Q/\alpha$. Alors le champ de vecteurs $\epsilon h_1
\frac{\partial}{\partial x}$ définit un automorphisme de $X_1$ qui est
trivial sur $X$ donné par $\phi_1(u)=u+\epsilon u' h_1$. 

On a alors 
\begin{multline*}
d\phi_1 \left(\delta_1(dy)\right)=d\phi_1(\alpha dx)=(\alpha+\epsilon
\alpha'h_1)(dx+\epsilon h'_1dx) =\\
\alpha dx+\epsilon(\alpha h_1)'dx= \delta_1(dx)+Qdx=\delta_2 dx
\end{multline*}
c'est-à-dire que $\delta_2=d\phi_1 \circ \delta_1$.

On peut donc remplacer $f_1$ par $f_1 \circ \phi_1$
et supposer que $df_1-df_2=\xi(\delta_1-\delta_2)=0$.

Comme $f$ est inséparable à la fibre spéciale, $\Aut_Y(Y_1)$ agit
trivialement sur $\delta_1$ et sur $df_1$. En effet, un élément de
$\Aut_Y(Y_1)$ est de la forme $Id+\epsilon \beta
\frac{\partial}{\partial y}$
et donc on a $$\delta_1(dy+\epsilon \beta dy)=f(1+\epsilon \beta)
\delta_1(dy)=(1+\epsilon \beta f'_0) \delta_1(dy)=\delta_1(dy)$$
car $f'_0=0$.

On a naturellement $f_1-f_2 \in (\epsilon)\otimes\Hom(\Omega_{Y_0/k},
f_* \O_{X_0})$. Décomposons à la fibre spéciale $f=s \circ i$ avec
$s:Z_0 \to Y_0$ étale et $i:X_0 \to Z_0$ purement inséparable. Comme
$df_1-df_2=0$ et que $i$ est de degré exactement $p$ on voit que
$f_1-f_2$ est en fait dans $(\epsilon)\otimes\Hom(\Omega_{Y_0/k}, s_*
\O_{Z_0})$ (on peut se ramener au cas des séries formelles par
localisation complétion, on voit alors que $f_1$ et $f_2$ coïncident
sur les monômes de degré premiers à $p$ car $df_1-df_2 = 0$, il
s'ensuit que $f_1-f_2$ est en fait une puissance $p$-ième,
c'est-à-dire dans l'image de $s_* \O_{Z_0}$).  Or $s$ est étale, donc
$\Hom(\Omega_{Y_0/k}, \O_{Y_0})$ agit transitivement sur 
$\Hom(\Omega_{Y_0/k}, s_* \O_{Z_0})$ (unicité de la déformation 
d'un morphisme étale).  Il existe donc un automorphisme de $Y$
envoyant $f_1$ sur $f_2$. 

On obtient alors le cas affine par recollement en utilisant le lemme
\ref{Automorphisme} et le fait que ce faisceau a un premier groupe de
cohomologie triviale. \qed

On remarquera que la démonstration précédente ne fonctionne pas si le
morphisme purement inséparable est de degré $p^2$ car, si $X=\Spec
k[[x]]$, $Y=\Spec k[[y]]$ et $f(y)=\sum a_i x^i$, alors l'information
différentielle donne les informations sur les $a_i$ avec $p \nmid i$,
la décomposition de $f$ en étale + séparable donne des informations
sur les $a_{p i}$ $p \nmid i$, mais on n'a pas d'information sur les
$a_{p^2 i}$ !

\begin{lemma}
Soit $A' \to A$ une petite extension et $\tilde \xi$ un relèvement de
$\xi$. Alors il existe un relèvement $(\tilde f, \tilde \delta, \tilde
\xi)$ de $(f, \delta, \xi)$ à $A'$.
\end{lemma}

\preuve Le morphisme $f$ se relève localement (car c'est un morphisme
localement d'intersection complète) et $\delta$ également (choix dune
base, ce qui est possible car $Z^1_{X/S}$ est localement libre). Il
s'agit alors de voir qu'on peut les choisir de sorte que
$df=a\delta$. 

Pour cela, regardons $df(dy)-a\delta(dy)=\epsilon d(Q)$ (ce qui est
possible car $df$ et $\delta$ sont tout deux exacts).
Par suite, changeant $f$ en $f+\epsilon Q \frac{\partial}{\partial
  x}$ (i.e. en composant avec l'automorphisme de $\tilde X$ définit
par $Q \frac{\partial}{\partial  x}$). Par suite, il existe un
relèvement au voisinage de chaque point. Comme le lieu des points $U$
ou $f_0$ se décompose en étale + inséparable est un ouvert  non vide,
ces relèvements locaux sont isomorphes en restriction à $U$ d'après le
lemme \ref{ClasseDeformation}. Ces relèvements locaux se laissent donc
recoller en une déformation globale. \qed

\begin{corollary}
\label{RelevementCaracteristitique0}
Soit $f\colon X \to Y$ un morphisme fini entre courbes propres et
lisses au-dessus d'un corps algébriquement clos $k$. Supposons que le
degré d'inséparabilité de $f$ est $\le p$ et que $f$ admette une
structure de morphisme rigidifié. Alors 
\begin{enumerate}
\item il existe un relèvement de $f$ sur $k[[t]]$ qui est
  génériquement séparable ;
\item il existe un anneau de valuation discrète $R$ d'inégales
  caractéristiques et de corps résiduel $k$, et un relèvement de $f$
  sur $R$.
\end{enumerate}
\end{corollary}

\preuve Le premier point est un corollaire immédiat du théorème
\ref{DeformationMorphismeRigidifie} (choisir $\xi \not = 0$).

Pour ce qui est du second, notons $g$ le genre de $X$ et $g'$ le genre
de $Y$ et considérons l'espace des modules $Mor_{g, g', \deg f}$
classifiant les morphismes de degré $\deg f$ entre courbes propres et
lisses de genre $g$ et $g'$. Cet espace est un champ algébrique
d'après la théorie des déformations des morphismes et le théorème
d'algébrisation d'Artin.

En particulier, $f$ définit un point de $Mor_{g, g', \deg
  f}$. D'après le premier point, il existe un morphisme  $\Spec k[[t]]
\to Mor_{g, g', \deg f}$ tel que l'image de $\Spec k((t))$ est dans
le lieu des morphismes séparables. D'après la théorie des déformations
des morphismes séparables entre courbes (qui est non obstruée), il
existe un anneau de valuation discrète $\tilde R$ d'inégales
caractéristiques et de corps résiduel $k((t))$, et un morphisme $\Spec
\tilde R \to Mor_{g, g', \deg f}$. Il s'ensuit que $f$ est la
spécialisation d'un point en caractéristique $0$, il existe donc un
anneau de valuation discrète $R$ dont le corps de fractions est une
extension finie de $\Frac \tilde R$ et dont le corps résiduel est $k$,
et un morphisme $\Spec R \to Mor_{g, g', \deg f}$ qui envoie le point
fermé sur $f$. \qed

La technique utilisée ici pour déformer en égales caractéristiques ne
permet pas de maîtriser $R$. On peut toutefois montrer, en utilisant
l'algébricité de $Mor_{g, g', \deg f}$, qu'on peut prendre pour
$R$ une extension finie de l'anneau des vecteurs de Witt $W(k)$. nous
conjecturons qu'il est possible de trouver un relèvement sur $W(k)$.

\section{Quelques conditions nécessaires}
\label{CN}

Nous souhaitons donner maintenant quelques conditions nécessaires pour
l'existence d'homomorphismes exacts : cf. les inéquations
\eqref{BorneRam} et \eqref{Bornea}). Celles ci sont assez grossières
mais nous semblent importantes pour montrer que le fait qu'un
morphisme puisse être muni d'une structure de morphisme rigidifié est
assez rare.

Donnons nous un morphisme séparable $h \colon Z \to Y$ entre courbes
propres et lisses au-dessus d'un corps $k$, notons $F\colon  X \to Z$
un morphisme purement inséparable de degré $p$ et supposons qu'il
existe un homomorphisme exact $\delta\colon (h\circ F)^*\Omega_{Y/k}
\to \Omega_{X/k}$. 

Notons de plus $r$ le degré du diviseur de ramification de $h$, $g$ le
genre de $X$ (et donc de $Z$), et $g'$ le genre de $Y$.
Comme $\delta$ est injectif, on a $2g-2-p\deg h(2g'-2) \ge 0$.

D'autre part, la formule de Hurwitz appliquée à $h$ nous donne
$$2g-2=\deg h(2g'-2)+r.$$

On obtient donc 

\begin{equation}
\label{BorneRam}
r \ge (2g-2)(1-\frac{1}{p}).
\end{equation}

En particulier, si $g$ est $\ge 2$, alors $h$ ne peut pas être étale.

Il est en fait possible d'être plus précis. En effet, comme vu dans le
lemme \ref{FaisceauHomfExct}, $\delta$ induit un homomorphisme
injectif $h^* \Omega_Y \to F_* B^1_X$.
On a donc une inclusion $H^0(Z, h^* \Omega_Y) \subset \H^0(X, B^1_X)$.

L'homomorphisme d'adjonction $\Omega_{Y/k} \to h_* h^*\Omega_{Y/k}$
permet alors de construire une injection  

$$\H^0(Y, \Omega_{Y/k}) \to \H^0(X, B^1_X).$$

La théorie des recouvrements infinitésimaux (cf. \cite{Milne},
Proposition 4.14) permet alors de voir que $\H^0(X,
B^1_X)=\H^1_{fl}(X, \aalpha_p)$. Notant $a$ la dimension de
$\H^1_{fl}(X, \aalpha_p)$, on obtient
\begin{equation}
\label{Bornea}
a \ge g'.
\end{equation}

En particulier, si la courbe $X$ est ordinaire (et donc également
$Z$), ce qui impose que $a=0$, et que $g' \ge 1$ alors $f$ ne peut
admettre de structure de morphisme rigidifié !

D'autre part, comme $\dim_k \H^0(Z, h^* \Omega_{Y/k}) \ge \deg h^*
\Omega_{Y/k}+1-g$ on trouve que 
\begin{equation}
a \ge g-1-r.
\end{equation}

\section{Exemple}
\label{CS}

Considérons une courbe propre et lisse $Y$ de genre $g'\ge 2$ qui est
supersingulière, choisissons un point $\p \in Y$ et considérons un 
$\ZZ/p^n\ZZ$-revêtement $Z \to Y$ étale en dehors de $\p$, et totalement
ramifié en $\p$ avec une différente $\dif$ qui peut être
arbitrairement grande (cf. \cite{KatzGabber} pour la construction de
tels revêtements) de sorte que le genre $g$ de $Z$ vérifie 
$$2g-2=p^n(2g'-2)+\dif$$
et $$1-p_Z=p^n(1-p_Y))-(p^n-1)$$
d'après \cite{CrewEtale} Corollary 1.8 (on note $p_Z$ et $p_Y$ les
$p$-rangs de $Z$ et de $Y$).

Comme $p_Y=0$ on trouve que $p_Z=0$, c'est-à-dire que $Z$ est
également supersingulière.

Notons $F\colon X \to Z$ le frobenius relatif. En particulier, $X$ est
également supersingulière d'après le théorème de descente des
morphismes étales par des morphismes purement inséparable. On a alors
\begin{multline*}
h^0(\Omega_{X} \otimes F^* h^* \Omega_Y^\vee) \ge 
\chi(\Omega_{X} \otimes F^* h^* \Omega_Y^\vee)= \\
\deg (\Omega_{X} \otimes F^* h^*
\Omega_Y^\vee)+1-g= \\
g-1-p^{n+1}(2g'-2)=(1-2p)p^n(g'-1)+\frac{\dif}{2}
\end{multline*}

En particulier, pour $\dif$ suffisamment grande ($p^n$ et $g'$ fixés), on
aura  $h^0(\Omega_{X} \otimes F^* h^* \Omega_Y^\vee) > 0$.

D'autre part, on a 
$$\H^0(X, \Omega_{X} \otimes F^* h^*
\Omega_Y^\vee)\overset{adjonction}{=}\H^0(Z, F_*(\Omega_{X}) \otimes h^* \Omega_Y^\vee)
\subset \H^0(Z, F_*(\Omega_{X}))=\H^0(X, \omega_X)$$
l'inclusion étant obtenue via le choix (non canonique) d'une injection
$\Omega_Y^\vee \to \O_Y$ qui existe car $g' \ge 2$.

Considérons maintenant l'opérateur de Cartier $C\colon \H^0(X,
\Omega_X) \to \H^0(X,\Omega_X)$. Comme $X$ est supersingulière,
l'opérateur $C$ est nilpotent. Par suite, comme $\H^0(X, \Omega_{X}
\otimes F^*h^* \Omega_Y^\vee) \subset \H^0(X, \Omega_{X})$
 est stable sous $C$ (l'opération de $C$ ne se fait que sur le premier
 terme de $\Omega_{X}\otimes F^*h^* \Omega_Y^\vee$),  il existe un
 élément $\delta$ non nul dans $\H^0(X, \Omega_{X} \otimes F^*h^*
 \Omega_Y^\vee)$ dont l'image par $C$ est nul, c'est-à-dire un élément
 de $\H^0(Z, F_*B^1_{X} \otimes h^* \Omega_Y^\vee)$.

\begin{remark}
Au moins si $n=1$, la procédure précédente permet de construire 
une famille non isotriviale de dimension $\lfloor\frac{\dif}{p}\rfloor$
(cf. la théorie d'Artin-Schreier). C'est-à-dire de dimension
$\frac{2g-2-p(2g'-2)}{p}$. 
\end{remark}

\bibliographystyle{../hamsalpha}

\providecommand{\bysame}{\leavevmode\hbox to3em{\hrulefill}\thinspace}

\end{document}